\documentclass[]{amsart}

\def\bi{\mathbb{I}}
\def\bj{\mathbb{J}}

\def\bn{\mathbb{N}}
\def\br{\mathbb{R}}
\def\bc{\mathbb{C}}

\def\h{\mathcal{H}}
\def\k{\mathcal{K}}
\def\l{\mathcal{L}}
\def\g{\mathcal{G}}
\DeclareMathOperator{\B}{\rm B}

\newcommand{\nor}[1]{#1_{\rm n}}

\newcommand{\po}[1]{#1^{\circ}}
\newcommand{\ppo}[1]{#1_{\circ}}
\newcommand{\du}[1]{#1^{\sharp}}
\newcommand{\pd}[1]{#1_{\sharp}}
\newcommand{\bdu}[1]{#1^{\sharp\sharp}}
\newcommand{\com}[1]{#1^{\prime}}
\newcommand{\ort}[1]{#1^{\perp}}

\newcommand{\norm}[1]{\Vert#1\Vert}
\newcommand{\cbnorm}[1]{\Vert#1\Vert_{{\rm cb}}}
\newcommand{\inner}[2]{\langle #1,#2\rangle}
\newcommand{\sumjn}{\sum_{j=1}^n}
\newcommand{\sumi}{\sum_{i\in\mathbb{I}}}
\newcommand{\sumj}{\sum_{j\in\mathbb{J}}}

\newcommand{\bb}[2]{{\rm B}(#1,#2)}
\newcommand{\abbb}[2]{{\rm B}_A(#1,#2)_B}
\newcommand{\bh}{{\rm B}(\mathcal{H})}
\newcommand{\bl}{{\rm B}(\mathcal{L})}

\newcommand{\ta}{\tilde{A}}

\newcommand{\bkh}{{\rm B}(\mathcal{K},\mathcal{H})}

\newcommand{\rown}[1]{{\rm R}_{n}(#1)}
\newcommand{\coln}[1]{{\rm C}_{n}(#1)}
\newcommand{\mat}[3]{{\rm M}_{#1,#2}(#3)}
\newcommand{\weakc}[1]{\overline{#1}}

\newcommand{\pro}{\stackrel{\gamma}{\otimes}}
\newcommand{\gc}{\gamma_{C}}
\newcommand{\gabc}[3]{\gamma^{#2}_{#1,#3}}

\newcommand{\npro}{\stackrel{\nu}{\otimes}_C}
\newcommand{\cnpro}{\stackrel{\nu}{\widetilde{\otimes}}}
\newcommand{\nc}{\nu_{C}}
\newcommand{\tnc}{\tilde{\nu}_{C}}
\newcommand{\Labc}[3]{\Lambda_{#1,#3}^{#2}}
\newcommand{\Lb}{\Lambda_B}
\newcommand{\labc}[3]{\lambda_{#1,#3}^{#2}}
\newcommand{\lb}{\lambda_B}

\newcommand{\arb}{{_A{\rm RM}_B}}
\newcommand{\ara}{{_A{\rm RM}_A}}
\newcommand{\anra}{{_A{\rm NRM}_A}}
\newcommand{\anrb}{{_A{\rm NRM}_B}}
\newcommand{\crc}{{{\rm CRM}_C}}
\newcommand{\coc}{{{\rm COM}_C}}
\newcommand{\cnrc}{{{\rm CNRM}_C}}

\DeclareMathOperator{\essup}{\rm essup}
\DeclareMathOperator{\eop}{{\rm ess}\oplus_{t\in\Delta}}

\DeclareMathOperator{\supp}{\rm supp}

\newtheorem{theorem}{Theorem}[section]
\newtheorem{lemma}[theorem]{Lemma}
\newtheorem{corollary}[theorem]{Corollary}
\newtheorem{proposition}[theorem]{Proposition}
\theoremstyle{remark}
\newtheorem{remark}[theorem]{Remark}
\theoremstyle{definition}
\newtheorem{definition}[theorem]{Definition}
\newtheorem{example}[theorem]{Example}
\numberwithin{equation}{section}

\begin{document}

\title[]{On tensor products of operator modules} 
\author{Bojan Magajna} 
\address{Department of Mathematics\\ University of Ljubljana\\
Jadranska 19\\ Ljubljana 1000\\ Slovenia}
\email{Bojan.Magajna@fmf.uni-lj.si}



\begin{abstract} 
The injective  tensor product of normal representable bimodules over von Neumann algebras
is shown to be  normal. The usual Banach module projective tensor product
of central representable bimodules over an Abelian C$^*$-algebra is shown 
to be representable. A
normal version of the projective tensor product is introduced  for  
central normal  bimodules.
\end{abstract}

\maketitle

\section{Introduction}

A Banach bimodule $X$ over C$^*$-algebras $A$ and $B$ is called {\em representable}
(\cite{AP}, \cite{Po}) if there exist Hilbert modules $\h$ and $\k$ over $A$ and $B$,
respectively (that is, Hilbert spaces with $*$-representations $\pi:A\to\bh$ and
$\sigma:B\to\B(\k)$) and an isometric $A,B$-bimodule homomorphism $X\to\bkh$. We
denote the class of all such bimodules by $\arb$, and by $\B_A(X,Y)_B$ the space
of all bounded $A,B$-bimodule maps from $X$ into $Y$. If, in addition, $A$ and
$B$ are von Neumann algebras and $\h$ and $\k$ are normal (that is, the 
representations $\pi$ and $\sigma$ are normal), then we say that $X$ is a 
{\em normal representable} $A,B$-bimodule, which we shall write as
$X\in\anrb$. 
In \cite{AP} the fundamentals of the tensor products of representable bimodules
are studied. In particular the {\em projective tensor seminorm} on the algebraic
tensor product $X\otimes_BY$ of two bimodules $X\in\arb$ and $Y\in{_B{\rm RM}_C}$
is defined by
\begin{multline}\label{11} \gabc{A}{B}{C}(w)=\inf\{\norm{\sumjn a_ja_j^*}^{1/2}
\norm{\sumjn b_j^*b_j}^{1/2}:\ w=\sumjn a_jx_j\otimes_By_jb_j,\\
a_j\in A,\ b_j\in B,\ x_j\in X,\ y_j\in Y,\ \norm{x_j}\leq1,\ \norm{y_j}\leq1
\}.
\end{multline}
Taking the quotient of
$X\otimes_BY$ by the zero space of this seminorm and completing, we obtain a
representable $A,C$-bimodule, denoted by $_AX\pro_BY_C$, and the induced norm 
on this bimodule is denoted by 
$\gabc{A}{B}{C}$ again. In the case $A=B=C=\bc$ this reduces to the usual projective tensor product
of Banach spaces, denoted simply by $X\pro Y$. As shown in \cite{AP}, this seminorm
can also be expressed by
\begin{equation}\label{12} \gabc{A}{B}{C}(\sumjn x_j\otimes_By_j)=\sup\sumjn
\theta(x_j,y_j),
\end{equation}
where the supremum is over all contractive bilinear maps $\theta$ from $X\times Y$
into $\B(\l,\h)$, with $\h$ and $\l$ cyclic Hilbert modules over $A$
and $C$ (respectively),  such that
$$\theta(axb,yc)=a\theta(x,by)c\ \ \mbox{for all}\ \ a\in A,\ b\in B,\ c\in C,\ 
x\in X,\ y\in Y.$$
Further \cite{AP}, the {\em injective tensor seminorm} is defined on $X\otimes_BY$
by
\begin{equation}\label{13} \Labc{A}{B}{C}(\sumjn x_j\otimes_By_j)=\sup\norm{\sumjn
\phi(x_j)\psi(y_j)},
\end{equation}
where the supremum is over all contractions $\phi\in\abbb{X}{\bkh}$ and
$\psi\in{\rm B}_B(Y,{\rm B}(\l,\k))_C$,  with $\h$, $\k$ and $\l$ 
cyclic Hilbert modules over $A$, $B$ and $C$ (respectively).  

\begin{remark}\label{acb} The  restriction that $\h$ and $\k$ in the 
above formulas are cyclic over $A$ and $B$ (respectively) implies by
an argument of Smith \cite[Theorem 2.1]{S} that each bounded  $A,B$-bimodule homomorphism $\phi$ 
from an operator 
$A,B$-bimodule into $\bkh$ is completely bounded with $\cbnorm{\phi}=\norm{\phi}$. 
Applying this to a pair $Y\subseteq X$ of representable $A,B$-bimodules 
and using the extension theorem
for completely bounded bimodule maps \cite{P}, \cite{W1}, it follows 
that each map $\phi\in\B_A(Y,\bkh)_B$ can be extended to a map $\psi\in\B_A(X,\bkh)_B$
with $\norm{\psi}=\norm{\phi}$. Thus in this respect such maps behave like linear
functionals.
\end{remark}

Clearly there are similar definitions of the `projective' and the 'injective' 
tensor seminorms (which turn out to be norms) in the category $\anrb$ for von
Neumann algebras $A$ and $B$; the only difference with the above definitions is that we
require the cyclic Hilbert modules $\h$, $\k$ and $\l$ to be normal. Now the natural
question is if these new norms are different from the above ones. In Section 2
we shall show that the two injective norms are  equal. Following the 
observation that the norm $\Labc{A}{B}{C}$ is in fact independent of $A$
and $C$, the proof of equality of the two
injective norms will be essentially a reduction to a density question concerning 
certain sets of normal states. Contrary to the injective, the two projective norms
are not the same even if $A=B=C$ is Abelian and the bimodules are central.
Here a $C$-bimodule $X$ is called {\em central} if $cx=xc$ for all $c\in C$ and $x\in X$.
We denote by $\crc$ the class of all central representable $C$-bimodules and
(if $C$ is a von Neumann algebra) by $\cnrc$ the subclass of all central normal
representable bimodules.

In Section 3 we show that ${_CX\pro_C Y_C}={_{\bc}X\pro_C Y_{\bc}}$ 
for all bimodules $X,Y\in\crc$. (Note that ${_{\bc}X\pro_C Y_{\bc}}$ is
just $X\pro_CY$, the quotient of the usual Banach space tensor product
$X\pro Y$ by the closed subspace generated
by all elements of the form $xc\otimes y-x\otimes cy$ ($x\in X,\ y\in Y,\ c\in C$)
\cite{R}.)
The main step of the proof will be to show that the central $C$-bimodule $X\pro_CY$ 
is representable, which in the more traditional terminology (see \cite{KiR})
means that the usual Banach space projective tensor product of $C$-locally convex 
modules over $C$ is already $C$-locally convex. This simplifies the corresponding
definition of such tensor product in \cite{KiR}. 

If $C$ is an Abelian von Neumann algebra and
$X,Y\in\cnrc$, the bimodule $Z=X\pro_CY$ is not necessarily normal. Therefore we  introduce in Section 4 a new tensor product
$X\npro Y$, which plays the role of the projective tensor product in the category
$\cnrc$. We show that $\nor{Z}:=X\npro Y$ is just the normal part of $Z$ in the sense
that each bounded $C$-bimodule map $\phi$ from $Z$ into a bimodule 
$V\in{_C{\rm NRM}_C}$
factors uniquely through $\nor{Z}$. Further, the norm of elements
in $\nor{Z}$ can be expressed by a formula similar to (\ref{11}), but involving 
infinite sums that
are not necessarily norm convergent. We do not know if there is an analogous 
formula in the case of non-central bimodules.

The background concerning operator spaces used implicitly in this 
article can be found in any of the books \cite{ER}, \cite{P}, \cite{Pi}.

\section{Normality of the injective operator bimodule tensor product}

If $A$, $B$ and $C$ are von Neumann algebras and $X\in\anrb$, 
$Y\in{_B{\rm NRM}_C}$, we  define a norm  on
$X\otimes_BY$ by 
\begin{equation}\label{21} \labc{A}{B}{C}(\sumjn x_j\otimes_By_j)=\sup\norm{\sumjn
\phi(x_j)\psi(y_j)},
\end{equation}
where the supremum is over all contractions $\phi\in\abbb{X}{\bkh}$ and
$\psi\in{\rm B}_B(Y,{\rm B}(\l,\k))_C$ with $\h$, $\k$ and $\l$ {\em normal} 
cyclic Hilbert modules over $A$, 
$B$ and $C$ (respectively). Except for the normality requirement on
Hilbert modules, this is the same formula as (\ref{13}), hence 
$\labc{A}{B}{C}\leq\Labc{A}{B}{C}$. 

\begin{remark} To show that $\labc{A}{B}{C}$ is definite, 
suppose that $w=\sumjn x_j\otimes y_j\in X\otimes_BY$ is such that
$\sumjn\phi(x_j)\psi(y_j)=0$ for all $\phi$ and $\psi$ as in the definition of
$\labc{A}{B}{C}$. We may assume that 
$X\subseteq\B(\h_B,\h_A)$
and $Y\subseteq\B(\h_C,\h_B)$ for some normal (faithful) Hilbert modules $\h_A$,
$\h_B$ and $\h_C$ over $A$, $B$ and $C$, respectively. Decomposing $\h_A$, 
$\h_B$ and $\h_C$ into direct sums of cyclic submodules,
$$\h_A=\oplus_i\com{e}_i\h_A,\ \ \ \h_B=\oplus_j\com{f}_j\h_B,\ \ \ \h_C=\oplus_k
\com{g}_k\h_C,$$
where $\com{e}_i\in\com{A}$, $\com{f}_j\in\com{B}$ and $\com{g}_k\in\com{C}$ are projections, and
considering the maps of the form  $\phi(x)=\com{e}_ix\com{b}\com{f}_j$ and 
$\psi(y)=\com{f}_j\com{b}y\com{g}_k$, where $\com{b}\in\com{B}$, it follows that
$[x_1,\ldots,x_n]\com{B}[y_1,\ldots,y_n]^T=0$, which implies that $\sumjn x_j
\otimes y_j=0$ (see e.g. \cite[Lemma 1.1]{M1}).
\end{remark}

We would like to show that $\Labc{A}{B}{C}=\labc{A}{B}{C}$,
but first we shall show that $\Labc{A}{B}{C}$ and $\labc{A}{B}{C}$
are  independent of $A$ and $C$.
We simplify the notation by $\lb:=\labc{\bc}{B}{\bc}$ and  $\Lb:=\Labc{\bc}{B}{\bc}$. 
Note that Remark \ref{acb}  implies
that both norms $\Labc{A}{B}{C}$ and $\labc{A}{B}{C}$  
are preserved under isometric embeddings of bimodules. 

The conjugate (that is, the dual) space $\h^*$ of a (left) Hilbert $A$-module $\h$ is 
regarded below as a right $A$-module by $\xi^*a=(a^*\xi)^*$ ($\xi\in\h,\ a\in A$),
where $\xi^*$ denotes $\xi$ regarded as an element of $\h^*$.

\begin{proposition}\label{p81} The seminorms $\Labc{A}{B}{C}$ (for representable
bimodules over C$^*$-algebras) and $\labc{A}{B}{C}$ (for normal representable
bimodules over von Neumann algebras)
do not depend on $A$ and $C$. 
\end{proposition}

\begin{proof}  Choose $\varepsilon>0$. Given 
$w=\sumjn x_j\otimes_By_j\in X\otimes_BY$ and contractions  
$\phi\in\abbb{X}{\bkh}$, $\psi\in{\rm B}_B(Y,{\rm B}(\l,\k))_C$ as in (\ref{13})
or (\ref{21}), 
we choose
unit vectors $\xi\in\h$ and $\eta\in\l$ such that
$$|\inner{\sumjn\phi(x_j)\psi(y_j)\eta}{\xi}|>\norm{\sumjn\phi(x_j)\psi(y_j)}-
\varepsilon.$$
Then 
$$\alpha:X\to\k^*,\ \alpha(x)=(\phi(x)^*\xi)^*\ \ \mbox{and}\ \ 
\beta:Y\to\k,\ \beta(y)=\psi(y)\eta$$
are contractive homomorphisms of modules over $B$  such that
$$|\sumjn\inner{\beta(y_j)}{\alpha(x_j)^*}|>\norm{\sumjn\phi(x_j)\psi(y_j)}-
\varepsilon.$$
This implies that $\Lb(w)\geq\Labc{A}{B}{C}(w)$ and $\lb(w)\geq\labc{A}{B}{C}(w)$.

To prove the inequality $\Lb(w)\leq\Labc{A}{B}{C}(w)$, let $\pi:B\to\B(\k)$ be a cyclic representation
and let $\alpha\in{\rm B}(X,\k^*)_B$, $\beta\in{\rm B}_B(Y,\k)$ be contractions such that
\begin{equation}\label{82} |\sumjn\inner{\beta(y_j)}{\alpha(x_j)^*}|>\Lb(w)-
\varepsilon.
\end{equation}
Since $\Labc{A}{B}{C}$ is preserved by inclusions we may assume that $X$ and $Y$ 
are  C$^*$-algebras
containing $A\cup B$ and $B\cup C$ (resp.). Then, since $\alpha$ and $\beta$
are complete contractions by a result of Smith quoted in Remark \ref{acb}, it 
follows by the representation theorem for such mappings
(see \cite[p. 102]{P}) that 
there exist Hilbert spaces $\h$ and $\l$,
$*$-representations $\Phi:X\to\bh$ and $\Psi:Y\to\bl$, unit vectors
$\xi\in\h$ and $\eta\in\l$ and contractions $S\in\bb{\k}{\h}$, $T\in\bb{\l}{\k}$
such that
\begin{equation}\label{83} \alpha(x)=\xi^*\Phi(x)S\ \ \mbox{and}\ \ 
\beta(y)=T\Psi(y)\eta.
\end{equation}
Clearly we may adjust $\h$, $\k$,  $S$ and $T$ so that $[\Phi(X)\xi]=\h$ and $[\Psi(Y)\eta]=\l$. Then it follows 
from (\ref{83})
(since $\alpha$ and $\beta$ are $B$-module maps) that
\begin{equation}\label{84} \Phi(b)S=S\pi(b)\ \ \mbox{and}\ \ 
T\Psi(b)=\pi(b)T\ \ (b\in B).
\end{equation}
Replace $\h$ with the subspace $\h_1=[\Phi(A)\xi]$ and $\l$ with
$\l_1=[\Psi(C)\eta]$ and define
$$\psi:Y\to\bb{\l_1}{\k},\ \ \ {\rm by}\ \ \ \psi(y)=T\Psi(y)|\l_1$$
and
$$\phi:X\to\bb{\k}{\h_1},\ \ \ {\rm by}\ \ \ \phi(x)=P\Phi(x)S,$$
where $P\in\bh$ is the orthogonal projection onto $\h_1$. Then $\eta\in\l_1$,
$\xi\in\h_1$ and by (\ref{83}) 
\begin{equation}\label{85} \alpha(x)=\xi^*\phi(x)\ (x\in X)\ \ \mbox{and}\ \ 
\beta(y)=\psi(y)\eta\ (y\in Y).
\end{equation}
Moreover, $\h_1$, $\k$ and $\l_1$ are cyclic  over $A$, $B$ and $C$ (respectively)  
and  (\ref{84}) (together with the fact that $\h_1$ and $\l_1$ are invariant under
$\Phi(A)$ and $\Psi(C)$, respectively)
implies that $\phi(axb)=\Phi(a)\phi(x)\pi(b)$ and 
$\psi(byc)=\pi(b)\psi(y)\Psi(c)$, thus 
$\phi\in\B_A(X,\B(\k,\h_1))_B$ and 
$\psi\in\B_B(X,\B(\l_1,\k))_C$ are of the type required in the definition of
the norm $\Labc{A}{B}{C}$.
Since from (\ref{85}) and (\ref{82}) we have that
$$\norm{\sumjn\phi(x_j)\psi(y_j)}\geq|\inner{\sumjn\phi(x_j)\psi(y_j)\eta}{\xi}|
=|\sumjn\inner{\beta(y_j)}{\alpha(x_j)^*}|>\Lb(w)-\varepsilon,$$
it follows that $\Labc{A}{B}{C}(w)\geq\Lb(w)$.

The proof of the inequality $\lb(w)\leq\labc{A}{B}{C}(w)$ is same as the proof
in the previous paragraph, with the addition that we must achieve that the
modules $\k$, $\h_1$ and $\l_1$ are normal. First, since $X\in\anrb$ and
$Y\in{_B{\rm NRM}_C}$ we may assume (by  standard  arguments)
that (up to isometric isomorphisms) $A,X,B\subseteq\B(\h_0)$ and $B,Y,C\subseteq\B(\l_0)$ for some Hilbert spaces
$\h_0$ and $\l_0$ (with the module multiplications just the products of operators).
Then (by Remark \ref{acb} again) we may assume that $X=\B(\h_0)$ and $Y=\B(\l_0)$.
By the definition of the norm $\lb$ we can choose a normal cyclic representation $\pi:B\to\k$ and $\alpha\in\B(X,\k^*)_B$,
$\beta\in\B_B(Y,\k)$ such that 
\begin{equation}\label{82n} |\sumjn\inner{\beta(y_j)}{\alpha(x_j)^*}|>
\lb(w)-\varepsilon.
\end{equation}
Let $U$ be the unit ball of $\B_B(Y,\k)={\rm CB}_B(Y,\k)$ (Remark \ref{acb}, 
note that $\k=\B(\bc,\k)$) and $U_{\sigma}$ the
weak* continuous maps in $U$. Since $Y=\B(\l_0)$, it follows from a variant of
\cite[2.5]{EK} that $U_{\sigma}$ is dense in $U$ in the point weak* topology; 
but since
$\k$ is reflexive, this topology has the same continuous linear functionals
as the point norm topology, 
hence by convexity $U_{\sigma}$ is dense in $U$ in
the point norm topology. With a similar result for $\B(X,\k^*)_B$, it follows that we may
assume that the maps $\alpha$ and $\beta$ in (\ref{82n}) are weak* continuous. 
But then the proof of the representation theorem for completely bounded mappings
\cite[Theorem 8.4]{P} (together with the Stinespring's construction) shows that the representations $\Phi$ and $\Psi$ constructed
in the previous paragraph are normal, which implies that the Hilbert modules
$\h_1$ and $\l_1$ over $A$ and $C$  are also normal. (Alternatively, 
we could just take the normal parts of $\Phi$ and $\Psi$...)
\end{proof}
Note that the analogy of Proposition \ref{p81} for the projective norm does
not hold, namely for a C$^*$-algebra $A$ the norm $\gabc{A}{\bc}{A}$ on 
$A\otimes A$ coincides with the Haagerup norm, while the norm 
$\gabc{\bc}{\bc}{\bc}$ is the usual Banach space projective tensor norm.

A subset $K$ of an $A,B$-bimodule $X$ is called {\em  $A,B$-absolutely convex}
if
$$\sumjn a_jx_jb_j\in K$$
for all $x_j\in K$ and $a_j\in A$, $b_j\in B$ satisfying $\sumjn a_ja_j^*\leq1$,
$\sumjn b_j^*b_j\leq1$.

\begin{lemma}\label{l82} If $K$ is a $B,\bc$-absolutely convex weak* compact
subset of a von Neumann algebra $B$, then the set $L=\{x^*x:\ x\in K\}$ is
convex and weak* compact.
\end{lemma}

\begin{proof} Given $x,y\in K$ and $t\in[0,1]$, consider the polar decomposition
$$\left[\begin{array}{c}
\sqrt{t}x\\
\sqrt{1-t}y
\end{array}\right]=\left[\begin{array}{c}
u\\
v
\end{array}\right]z,$$
where $z=\sqrt{tx^*x+(1-t)y^*y}$ and $[u,v]^T$ is the partial isometric part.
Since 
$$z=[u^*\ v^*]\left[\begin{array}{c}
\sqrt{t}x\\
\sqrt{1-t}y
\end{array}\right]=u^*x\sqrt{t}+v^*y\sqrt{1-t}$$
and $K$ is $B,\bc$-absolutely convex, $z\in K$. It follows that
$tx^*x+(1-t)y^*y=z^*z\in L$, proving that $L$ is convex.

Since $K$ (hence also $L$) is bounded, it suffices now to prove that $L$ is
closed in the strong operator topology (SOT). Let $y$ be in the closure of $L$ 
and $(x_j)$ a net in $K$ such that $(x_j^*x_j)$ converges to $y$ in the
SOT. Since the function $x\mapsto\sqrt{x}$ is SOT continuous on bounded
subsets of $B^+$, the net $(|x_j|)$ converges to $\sqrt{y}$. Since $K$
is $B,\bc$-absolutely convex, the polar decomposition shows that $|x_j|\in K$.
Since $K$ is weak* closed, it follows that $\sqrt{y}\in K$, hence $y\in L$.
\end{proof}

We denote by $\rown{B}$ and $\coln{B}$ the set of all $1\times n$ and
$n\times 1$ matrices (respectively) with the entries in a set $B$.
\begin{theorem}\label{t83} For all $X\in{\rm NRM}_B$ and $Y\in{_B{\rm NRM}}$,
$\Lb=\lb$ on $X\otimes_BY$.
\end{theorem}

\begin{proof} The theorem will be proved first for free modules by translating
the problem to states on $B$ and approximating states by normal states. Then
elements of general modules will be approximated by elements of free modules.

First assume that $X$ and $Y$ are free with basis $\{x_1,\ldots,
x_n\}$ and $\{y_1,\ldots,y_n\}$, respectively. More precisely, set
$$x=[x_1,\ldots,x_n],\ \ \ y=[y_1,\ldots,y_n]^T$$
and assume that the two maps
$$f:\coln{B}\to X,\ f(b)=xb\ \ \mbox{and}\ \ g:\rown{B}\to Y,\ g(b)=by$$
are invertible (with bounded inverses by the open mapping theorem). Set
$${\mathcal{S}}=\{b\in\coln{B}:\ \norm{xb}\leq1\},\ \  \mathcal{T}=\{
b\in\rown{B}:\ \norm{by}\leq1\}$$
and 
$$\alpha=\sup\{\norm{b}:\ b\in\mathcal{S}\cup\mathcal{T}\}.$$
Let $0<\varepsilon<1$. Choose $w\in X\otimes_BY$ and note that $w$ can be written as 
\begin{equation}\label{4w} w=\sum_{i,j=1}^nx_i\otimes_Bd_{ij}y_j\ \ 
(d_{ij}\in B).
\end{equation}
By the definition of $\Lb$ there exist a cyclic representation
$\pi:B\to\k$ and contractions $\phi\in{\rm B}(X,\k^*)_B$, $\psi\in
{\rm B}_B(Y,\k)$ such that
\begin{equation}\label{86} |\sum_{i,j=1}^n\inner{\pi(d_{ij})\psi(y_{j})}
{\phi(x_i)^*}|>\Lb(w)-\varepsilon.
\end{equation}
Let $\xi_0\in\k$ be a unit cyclic vector for $\pi(B)$, $\rho$ the state
$\rho(b)=\inner{\pi(b)\xi_0}{\xi_0}$ on $B$, and choose $a_i,c_i\in B$ so that
\begin{equation}\label{87} \norm{\phi(x_i)^*-\pi(a_i^*)\xi_0}<\varepsilon\ \ 
\mbox{and}\ \ \norm{\psi(y_i)-\pi(c_i)\xi_0}<\varepsilon\ \ (i=1,\ldots,n).
\end{equation}
For $b=[b_{ij}]\in\mat{m}{n}{B}$ denote the matrix
$[\pi(b_{ij})]$ simply by $\pi(b)$.
Set
\begin{equation}\label{88} \xi=[\phi(x_1)^*,\ldots,\phi(x_n)^*]^T\ \ (\in\k^n),\ \ 
\ \eta=[\psi(y_1),\ldots,\psi(y_n)]^T\ \ (\in\k^n),
\end{equation}
$$a=[a_1,\ldots,a_n]\ \ \mbox{and}\ \ c=[c_1,\ldots,c_n]^T.$$
Then from (\ref{87})
\begin{equation}\label{89} \norm{\xi-\pi(a)^*\xi_0}<\varepsilon\sqrt{n}\ \
\mbox{and}\ \ \norm{\eta-\pi(c)\xi_0}<\varepsilon\sqrt{n}.
\end{equation}
Since  $\psi$ is a contractive
$B$-module map, we have $$\norm{\sumjn\pi(b_j)\eta_j}=\norm{\sumjn\pi(b_j)\psi(y_j)}
=\norm{\psi(\sumjn b_jy_j)}\leq\norm{\sumjn b_jy_j},$$ hence (and similarly)
\begin{equation}\label{810} \norm{\pi(b)^*\xi}\leq\norm{xb}\ (b\in\coln{B})\ \ \ 
\mbox{and}\ \ \ \norm{\pi(b)\eta}\leq\norm{by}\ (b\in\rown{B}).
\end{equation}
Thus, if $b\in\mathcal{S}$, then
$$\begin{array}{lll}
\rho(abb^*a^*)&=&\norm{\pi(b^*a^*)\xi_0}^2\\
 &\leq&\left(\norm{\pi(b)^*\xi}+\norm{\pi(b)^*(\pi(a)^*\xi_0-\xi)}\right)^2\\
 &\leq&\left(\norm{xb}+\norm{\pi(b)}\varepsilon\sqrt{n}\right)^2\ \mbox{(by
 (\ref{810}) and (\ref{89})}\\
 &\leq&(1+\alpha\varepsilon\sqrt{n})^2\ \mbox{(by definition of}\ \mathcal{S}
 \ \mbox{and}\ \alpha)\\
 &=&:\beta.
 \end{array}$$
Similar arguments are valid for $b\in\mathcal{T}$, hence
 \begin{equation}\label{811} \rho(abb^*a^*)\leq\beta\ (b\in\mathcal{S})\ \ \ 
 \mbox{and}\ \ \ \rho(c^*b^*bc)\leq\beta\ (b\in\mathcal{T}).
 \end{equation}
 Set
 $$K_1=\{b^*a^*:\ b\in\mathcal{S}\},\ \ K_2=\{bc:\ b\in\mathcal{T}\},$$ 
 $$L_1=\{v^*v:\ v\in K_1\},\ \ L_2=\{v^*v:\ v\in K_2\}.$$
 Since $X$ and $Y$ are normal modules over $B$, $\mathcal{S}$ and $\mathcal{T}$
 are weak* closed; moreover, since $f$ and $g$ are invertible, $\mathcal{S}$ and 
 $\mathcal{T}$ are bounded, hence
 weak* compact. Thus,  $K_1$ and
 $K_2$ are also weak* compact. To verify that the subset $\mathcal{T}$ of 
 $\rown{B}$ is  $B,\bc$-absolutely
 convex, let $b_j\in\mathcal{T}$ ($j=1,\ldots,n$) and
 let $\lambda_j\in\bc$ and $d_j\in B$ satisfy $\sum|\lambda_j|^2\leq1$ and 
 $\sum d_jd_j^*\leq1$. Then to show that
 $\sum(d_jb_j\lambda_j)$ is in $\mathcal{T}$, just note
 that 
 $\norm{(\sum d_jb_j\lambda_j)y}=\norm{\sum d_j(b_jy)\lambda_j}\leq\max_j
 \norm{b_jy}\leq1$. Similarly $\mathcal{S}$ is $\bc,B$-absolutely convex
 and it follows that $K_1$ and $K_2$ are
 $B,\bc$-absolutely convex.
 
Now we deduce by Lemma \ref{l82} that $L_1$ and $L_2$ are convex weak*
compact subsets of $B_h$ (the self-adjoint part of $B$), hence the same holds for the convex
hull ${\rm co}(L_1\cup L_2)$ and therefore the set
$$L={\rm co}(L_1\cup L_2)-B^+$$
is weak* closed since $B^+$ (the positive part of $B$) is weak* closed. Set
$$\po{L}=\{\theta\in\du{B}:\ {\rm Re}(\theta(v))\leq1\ \forall{v}\in L\}\ \ 
\mbox{and}\ \ \ppo{L}=\po{L}\cap\pd{B}.$$
Since $L$ is weak* closed and convex, $\ppo{L}$ is weak* dense in $\po{L}$ by a 
variant of the bipolar theorem. From (\ref{811}) we have that $\rho\in\beta(\po{L}_1
\cap\po{L}_2)=\beta\po{({\rm co}(L_1\cup L_2))}$, hence (since $\rho$ is positive)
$\rho\in\beta\po{L}$. Since $\ppo{L}$ is weak* dense in
$\po{L}$, there exists an $\omega_0\in\beta\ppo{L}$ such that 
\begin{equation}\label{812} |(\omega_0-\rho)(\sum_{i,j=1}^na_id_{ij}c_j)|
<\varepsilon\ \ \mbox{and}\ \ |(\omega_0-\rho)(1)|<\varepsilon.
\end{equation}
(Here $d_{ij}$ are as in (\ref{4w}), thus $d_{ij}$, $a_i$ and $c_j$ are fixed.)
Since $L\supseteq-B^+$ and $\omega_0\in\beta\ppo{L}$, $\omega_0$ is positive, hence $\omega=\omega_0/
\omega_0(1)$ is a state. Since $\norm{\omega-\omega_0}=\norm{(1-\omega_0(1))\omega}
=|1-\omega_0(1)|<\varepsilon$, we have from (\ref{812}) that
\begin{equation}\label{813} |(\omega-\rho)(\sum_{i,j=1}^na_id_{ij}c_j)|<
D\varepsilon,
\end{equation}
where $D=1+\norm{\sum_{i,j=1}^na_id_{ij}c_j}$.
Let $\sigma:B\to\bh$ be the normal representation constructed from $\omega$
by the GNS construction and let $\eta_0\in\h$ be the corresponding unit cyclic
vector. From 
(\ref{813}), (\ref{87}) and (\ref{86}) we deduce that
\begin{equation}\label{814}\begin{array}{lll}|\sum_{i,j=1}^n\inner{\sigma(a_id_{ij}c_j)\eta_0}
{\eta_0}|&=&|\omega(\sum_{i,j=1}^na_id_{ij}c_j)|\\
 &>&|\rho(\sum_{i,j=1}^na_id_{ij}c_j)|-D\varepsilon\\
 &=&|\sum_{i,j=1}^n\inner{\pi(a_id_{ij}c_j)\xi_0}{\xi_0}|-D\varepsilon\\
&>&|\sum_{i,j=1}^n\inner{\pi(d_{ij})\psi(y_j)}{\phi(x_i)^*}|-D\varepsilon\\
& & -
n^2\varepsilon\max_{i,j}\norm{d_{ij}}(\norm{x}+\norm{y}+\varepsilon)\\
&>&\Lb(w)-r(\varepsilon),\end{array}
\end{equation}
where $r(\varepsilon)$ tends to $0$ as $\varepsilon\to 0$.

Define $\Phi_0\in{\rm B}(X,\h^*)_B$ and $\Psi_0\in{\rm B}_B(Y,\h)$ by
\begin{equation}\label{815} \Phi_0(\sumjn x_jb_j)=(\sumjn\sigma(b_j^*a_j^*)
\eta_0)^*,\  \  \Psi_0(\sumjn b_jy_j)=\sumjn\sigma(b_jc_j)\eta_0\
\ (b_j\in B).
\end{equation}
Since $\omega_0\in\beta\ppo{L}$,  $\omega=\omega_0/\omega_0(1)$ and $\norm{\omega-\omega_0}<
\varepsilon$, we have that $\omega\in\omega_0(1)^{-1}\beta\ppo{L}\subseteq
(1-\varepsilon)^{-1}\beta\ppo{L}$, hence it follows from (\ref{815}) 
(noting that $abb^*a^*\in L$ if $b\in\mathcal{S}\subseteq\coln{B}$) that
$$\norm{\Phi_0(xb)}^2=\norm{\sigma(b^*a^*)\eta_0}^2=\omega(abb^*a^*)\leq(1-
\varepsilon)^{-1}\beta\ \ (b\in\mathcal{S})$$
and similarly
$$\norm{\Psi_0(by)}^2\leq(1-\varepsilon)^{-1}\beta\ \ (b\in\mathcal{T}).$$
Thus, with $\delta=(1-\varepsilon)^{-1/2}\beta^{1/2}=
(1-\varepsilon)^{-1/2}(1+\alpha\varepsilon\sqrt{n})$, we have (recalling the definitions of $\mathcal{S}$
and $\mathcal{T}$) that $\norm{\Phi_0}\leq\delta$ and
$\norm{\Psi_0}\leq\delta$. From (\ref{815}), $\Phi_0(x_j)=
(\sigma(a_j^*)\eta_0)^*$
and $\Psi_0(y_j)=\sigma(c_j)\eta_0$, hence we may rewrite (\ref{814}) as
$$|\sum_{i,j=1}^n\inner{\sigma(d_{ij})\Psi_0(y_j)}{\Phi_0(x_i)^*}|>\Lb(w)-r(\varepsilon).$$
Finally, setting $\Phi=\frac{1}{\delta}\Phi_0$ and $\Psi=\frac{1}{\delta}\Psi_0$,
we have  a normal cyclic Hilbert module $\h$ and contractions
$\Phi\in{\rm B}(X,\h^*)_B$, $\Psi\in{\rm B}_B(Y,\h)$ such that
$|\sum\inner{\sigma(d_{ij})\Psi(y_j)}{\Phi(x_i)^*}|$ approaches 
$\Lb(w)$ as $\varepsilon$ tends to $0$ since
$r(\varepsilon)\to 0$ and $\delta\to1$. Thus $\Lb(w)=\lb(w)$.

In general, when $X$ and $Y$ are not free, let
$w=\sumjn x_j\otimes_By_j\in X\otimes_BY$ and 
$$X_1=X\oplus\rown{B}\ \ \mbox{and}\ \ Y_1=Y\oplus\coln{B}.$$
Since both norms $\Lb$ and $\lb$ respect isometric embeddings, it suffices
to prove that $\Lb(w)\leq\lb(w)$ in $X_1\otimes_BY_1$. For each real $t>0$ put
$$w(t)=\sumjn(x_j,te_j^T)\otimes_B(y_j,te_j),$$
where $e_j=(0,\ldots,1,\ldots,0)\in\coln{\bc}\subseteq\coln{B}$. Since the elements
$x_j(t):=(x_j,te_j^T)$ ($j=1,\ldots,n$) generate a free module in the
above sense and similarly the $y_j(t):=(y_j,te_j)$, it follows 
that $\Lb(w(t))=\lb(w(t))$.
But, as $t$ tends to $0$, $\Lb(w(t))$ tends to $\Lb(w)$ 
(since $\Lb(w(t)-w)\leq t\sumjn(\norm{x_j}+\norm{y_j}+t)$) and 
$\lb(w(t))$ tends to $\lb(w)$, hence $\Lb(w)=\lb(w)$.
\end{proof}

By Theorem \ref{t83} and Proposition \ref{p81} the injective norm is given by
(\ref{21}) where $\h$, $\k$ and $\l$ are normal, hence using the condition
for normality recalled in the last part of Theorem \ref{t59} below
we conclude:

\begin{corollary} If $X\in\anrb$ and $Y\in{_B{\rm NRM}_C}$, then 
$X\stackrel{\Lambda}{\otimes}_BY\in{_A{\rm NRM}_C}$.
\end{corollary}

\section{The projective tensor product of central bimodules}

{\em Throughout  this section $C$ is a unital Abelian C$^*$-algebra, 
$\tilde{C}$ the universal von Neumann envelope of $C$ in the standard
form and
$X,Y\in\crc$.}

\begin{remark}\label{r47} For an Abelian C$^*$-algebra $C$ we 
 denote by $\Delta$ the spectrum of $C$ and
 by $C_t$ the kernel of a character $t\in \Delta$. For a bimodule $X\in\crc$ we 
 consider the quotients 
 $X(t)=X/[C_tX]$. Given $x\in X$ we denote by $x(t)$ the coset 
 of $x$ in $X(t)$. It is known (see \cite[p. 37, 41]{DG} and \cite[p.71]{Po}
 or \cite{M10}) that
 the function
 \begin{equation}\label{47}\Delta\ni t\mapsto\norm{x(t)}\end{equation}
 is upper semicontinuous and that
 \begin{equation}\label{471}\norm{x}=\sup_{t\in\Delta}\norm{x(t)}.\end{equation}
 We shall call the embedding $$X\to\oplus_{t\in\Delta} X(t),\ \ 
 x\mapsto(x(t))_{t\in\Delta}$$
 the {\em canonical decomposition} of $X$. 
 \end{remark}
 
Let $X\pro_C Y$ be the quotient of the Banach space projective tensor product
$X\pro Y$ by the closed subspace generated
by all elements of the form $xc\otimes y-x\otimes cy$ ($x\in X$, $y\in Y$,
$c\in C$). First we shall prove  that $X\pro_CY$
is a representable $C$-bimodule. In classical terminology, this means that
$X\pro_CY$ is $C$-locally convex, which simplifies the definition of 
the tensor product of $C$-locally convex modules  \cite[p. 445]{KiR} since it
eliminates the need for Banach bundles.

Consider the canonical decompositions $X\to\oplus_{t\in\Delta}X(t)$
and $Y\to\oplus_{t\in\Delta}Y(t)$ along the spectrum $\Delta$ of $C$ (see
Remark \ref{r47}). For
each $t\in\Delta$ the $C$-balanced bilinear map
$$\kappa_t:X\times Y\to X(t)\pro Y(t),\ \ \kappa_t(x,y)=x(t)\otimes y(t)$$
induces a contraction 
$\tilde{\kappa}_t:X\pro_CY\to X(t)\pro Y(t)$. Since the kernel of $\tilde{\kappa}_t$
contains the submodule $C_t(X\pro_CY)$ (where $C_t=\ker{t}$), 
$\tilde{\kappa}_t$ induces a contraction
$$ \mu_t:(X\pro_CY)(t)\to X(t)\pro Y(t).$$
On the other hand, the natural bilinear map $X\times Y\to X\otimes Y\to
(X\pro_CY)(t)$ annihilates
$C_tX\times Y$ and $X\times C_tY$, hence it induces a bilinear map
$X(t)\times Y(t)\to(X\pro_CY)(t)$ and therefore a linear map
$\sigma_t:X(t)\pro Y(t)\to (X\pro_CY)(t)$, which must be a contraction by
the maximality of the cross norm $\gamma$. Clearly $\sigma_t$ is inverse
to $\mu_t$ and since both are contractions, they must be isometries. Thus,
we have the isometric identification
\begin{equation}\label{72} (X\pro_CY)(t)=X(t)\pro Y(t)\ \ (t\in\Delta).
\end{equation}
For each $w\in X\pro_CY$ we denote by $w(t)$ the corresponding class in
$X(t)\pro Y(t)$.
 
We begin with the following result. 
\begin{theorem}\label{t71} The natural contraction 
\begin{equation}\label{73}\kappa:X\pro_CY\to\oplus_{t\in\Delta}(X(t)\pro Y(t)),
\ \ \kappa(x\otimes_Cy)=
(x(t)\otimes y(t))_{t\in\Delta}
\end{equation}
is isometric, hence $X\pro_CY$ is a representable $C$-bimodule.
\end{theorem}
For the proof we need some preparation. Set 
$Z=X\pro_CY$. Since the $C$-bimodule $\oplus_{t\in\Delta}Z(t)$ is clearly 
representable and $Z(t)=X(t)\pro Y(t)$ by (\ref{72}), it will suffice to prove
that the map (\ref{73}) is isometric.
Further, since for each element $w\in X\pro_CY$ its norm is equal to
$$\gc(w)=\sup\{|\theta(w)|:\ \theta\in\du{(X\pro_CY)},\ \norm{\theta}\leq1\},$$
it will suffice to show that
\begin{equation}\label{74} |\theta(w)|\leq\sup_{t\in\Delta}\norm{w(t)}
\end{equation}
for each $\theta$ in the unit ball of $X\pro_CY$.

\begin{remark}[Definition] Given $\theta\in\du{(X\pro_CY)}$ (regarded as a bilinear form) and an open subset $\Lambda$ of $\Delta$, let us define
that 
$$\theta|\Lambda=0\Longleftrightarrow \theta(x,cy)=0\ \forall c\in C=C(\Delta)\ 
\mbox{with}\ {\rm supp}\, c\subseteq\Lambda\ \mbox{and}\ \forall x\in X, \forall y\in Y.$$
If $(\Lambda_j)$ is a family of open subsets of $\Delta$ with the union $\Lambda$
and if $\theta|\Lambda_j=0$ for all $j$, then a standard partition of unity argument
shows that $\theta|\Lambda=0$. It follows that there exists the largest open
subset $\Lambda$ of $\Delta$ such that $\theta|\Lambda=0$; then $\Delta\setminus
\Lambda$ is called the {\em support of $\theta$}, denoted by $\supp\theta$.
\end{remark}

\begin{lemma}\label{l72} If $\theta$ is an extreme point of the unit ball of 
$\du{(X\pro_CY)}$ then $\supp\theta$ is a singleton.
\end{lemma}

\begin{proof} We can extend $\theta$ to a contractive bilinear form on 
$\bdu{X}\times\bdu{Y}$, denoted by $\theta$ again, such that the maps
\begin{equation}\label{720} \bdu{X}\ni F\mapsto\theta(F,y)\ (y\in Y)\ \ \mbox{and}\ \  
\bdu{Y}\ni G\mapsto\theta(x,G)\ (x\in X)
\end{equation}
are weak* continuous (see \cite[p. 12]{DF} if necessary). Since $X$ and $Y$
are representable,  we may
regard $\bdu{X}$ and $\bdu{Y}$ as normal dual  bimodules over
$\tilde{C}=\bdu{C}$ by \cite{M10} (this is  explained in more detail
also in the beginning of Section 4). 
In particular, for each bounded Borel function $f$ on $\Delta$ and each $y\in Y$,
$fy$ is defined as an element of $\bdu{Y}$. Thus, we may define a
bilinear form $f\theta$ on $X\times Y$ by
$$(f\theta)(x,y)=\theta(x,fy),$$
which satisfies  
\begin{equation}\label{75} (cf)\theta=c(f\theta)\ \ (c\in C).
\end{equation}
Using the separate weak* continuity of the maps (\ref{720}) and the fact that
the $\tilde{C}$-bimodules $\bdu{X}$ and $\bdu{Y}$ are normal, it also follows that
\begin{equation}\label{730} \theta(xc,y)=\theta(x,cy)\ \ (c\in\tilde{C},\ 
x\in X,\ y\in Y).
\end{equation}

Suppose that there exist two different points
$t_1,t_2\in\supp\theta$. Choose an open neighborhood $\Delta_1$ 
of $t_1$ such that $t_2\notin\overline{\Delta}_1$ and let $\chi$ be the 
characteristic function of $\Delta_1$. Then
$\chi\theta\ne0$.  (Indeed, $\chi\theta=0$ would imply
for all $c\in C$ with support in $\Delta_1$ that $c\theta=(c\chi)\theta=
c(\chi\theta)=0$ by (\ref{75}), hence $\theta(x,cy)=
(c\theta)(x,y)=0$ for all $x,y$, thus $\theta|\Delta_1=0$, but this would
contradict the fact that $t_1\in{\rm supp}\,\theta$.) 
Similarly $(1-\chi)\theta\ne0$.
Further, 
\begin{equation}\label{76} \norm{\chi\theta}+\norm{(1-\chi)\theta}=\norm{\theta}=1.
\end{equation}
Indeed, given $x,u\in X$ and $y,v\in Y$, for suitable $\alpha,\beta\in\bc$ of
modules $1$ we compute by using the property (\ref{730}) that
$$\begin{array}{lll}
|(\chi\theta)(x,y)|+|((1-\chi)\theta)(u,v)|&=&\alpha(\chi\theta)(x,y)+\beta((1-\chi)  
\theta)(u,v)\\
&=&\theta(x\chi,\alpha\chi y)+\theta(u(1-\chi),\beta(1-\chi)v)\\
&=&\theta(x\chi+u(1-\chi),\alpha\chi y+\beta(1-\chi)v)\\
&\leq&\norm{x\chi+u(1-\chi)}\norm{\alpha\chi y+\beta(1-\chi)v}\\
&\leq&\max\{\norm{x},\norm{u}\}\max\{\norm{y},\norm{v}\}.
 \end{array}$$
This implies that $\norm{\chi\theta}+\norm{(1-\chi)\theta}\leq1\ (=\norm{\theta})$,
while the reverse inequality is immediate from $\theta=\chi\theta+(1-\chi)\theta$.

Setting $s=\norm{\chi\theta}$, it follows  that $\theta$ is the
convex combination 
$\theta=s(s^{-1}\chi\theta)+(1-s)((1-s)^{-1}(1-\chi)\theta),$
where $s^{-1}\chi\theta$ and (by (\ref{76})) $(1-s)^{-1}(1-\chi)\theta$ are in the unit ball
of $\du{(X\pro_CY)}$. This is a contradiction since $\theta$ is an extreme point.
\end{proof}

\begin{proof}[Proof of Theorem \ref{t71}.] As we have already noted, it suffices 
to prove (\ref{74}). By the Krein Milman theorem we may assume that $\theta$ is an extreme
point in the unit ball of $X\pro_CY$. Then by Lemma \ref{l72} $\supp\theta=
\{t\}$ for some $t\in\Delta$. This implies that $\theta(XC_t,Y)=0=\theta(X,C_tY)$
since each $c\in C_t$  can be approximated by functions with 
supports in $\Delta\setminus
\{t\}$. Consequently $\theta$
can be factored through $X(t)\times Y(t)$, in other words, there exists a contraction 
$\theta_t\in\du{(X(t)\pro Y(t))}$ such that $\theta=\theta_t\circ\tilde{\kappa}_t$. 
It follows
that $|\theta(w)|\leq\norm{w(t)}$ for each $w\in X\pro_CY$.
\end{proof}

\begin{remark}\label{rmd} If $Z\in\crc$, then $\norm{w}=\sup\{\norm{\phi(w)}:\ 
\phi\in\B_C(Z,\tilde{C}),\ \norm{\phi}\leq1\}$ (this is known, \cite{Po});
moreover, if $Z\in\cnrc$, then we may replace in this formula $\tilde{C}$
by $C$. The later fact can be deduced from \cite{M10} by identifying the
proper bimodule dual of $Z$ with $\B_C(Z,C)$, but can also be deduced from
an earlier result of Halpern \cite[Theorem 3]{Ha} by representing $Z$ (and $C$) in some
$\B(\h)$ and noting that then $Z\subseteq\com{C}$ since $Z$ is central.
\end{remark}

\begin{corollary}\label{c73} For
each $w\in X\otimes_CY$
\begin{equation}\label{77}\gc(w)=\inf\{\norm{\sumjn c_j}:\ 
w=\sumjn c_jx_j\otimes_Cy_j,\ c_j\in C^{+},\ x_j\in B_X,\ y_j\in B_Y\},
\end{equation}
hence ${_{\bc}X\pro_CY_{\bc}}={_{C}X\pro_CY_{C}}$ and this is just the usual 
projective tensor product $X\pro_CY$ of Banach $C$-modules.
\end{corollary}

\begin{proof} Since by Theorem \ref{t71} $X\pro_CY\in\crc$, by Remark \ref{rmd}
the norm of $w\in X\pro_CY$ is
$\gc(w)=\sup\{\norm{\phi(w)}:\ \phi\in{\rm B}_C(X\pro_CY,\tilde{C}),\ \norm{\phi}
\leq1\}.$
For $w$ of the form $w=\sumjn c_jx_j\otimes_Cy_j$,
where $c_j\in C^+$, $\norm{x_j}\leq1$, $\norm{y_j}\leq1$, and a contraction
$\phi\in{\rm B}_C(X\pro_CY,\tilde{C})$ we have
$$\begin{array}{l} \norm{\phi(w)}=\norm{\sumjn c_j^{1/2}\phi(x_j\otimes_Cy_j)c_j^{1/2}}\\
\leq\|[c_1^{1/2},\ldots,c_n^{1/2}]\|\max_j\norm{\phi(x_j\otimes_Cy_j)}
\|\left[\begin{array}{l}c_1^{1/2}\\ \vdots\\ c_n^{1/2}\end{array}\right]\|\\
 \leq\norm{\sumjn c_j}\max_j\norm{x_j\otimes_Cy_j}\leq\norm{\sumjn c_j}.
\end{array}$$
This implies that $\gc(w)$ is dominated by the right side of (\ref{77}). But,
by definition
$$\gc(w)=\inf\{\sumjn\lambda_j:\ w=\sumjn\lambda_jx_j\otimes_Cy_j,\ \lambda_j\in
\br^+,\ x_j\in B_X,\ y_j\in B_Y,\ n\in\bn\},$$
which clearly dominates the right side of (\ref{77}) since $\bc\subseteq C$. The
conclusions of the corollary follow now from definitions of the 
corresponding norms.
\end{proof}

\begin{example}\label{e76} If $C$ is an Abelian von Neumann algebra and
$X,Y\in\cnrc$, then the representable $C$-bimodule $X\pro_CY$ is not necessarily  
normal. 
To show this, we modify an idea from \cite[Example 3.1]{KW}. Let $U_0\subseteq U$ and
$V$ be Banach spaces such that the  contraction $U_0\pro V\to
U\pro V$ is not isometric. Choose
$t_0\in\Delta$ and set 
$X=\{f\in C(\Delta,U):\ f(t_0)\in U_0\},$ $Y=C(\Delta,V).$
Then
$$X(t)=\left\{\begin{array}{ll}
U&\mbox{if}\ t\ne t_0\\
U_0&\mbox{if}\ t=t_0
\end{array}\right.
\ \ \ \mbox{and}\ \ \ Y(t)=V\ \mbox{for all}\ t\in\Delta.$$
Choose $w=\sumjn u_j\otimes v_j\in U_0\otimes V$ so that $\norm{w}_{U\pro V}
<\norm{w}_{U_0\pro V}$, denote by $\tilde{u}_j$ and $\tilde{v}_j$ the
constant functions $\tilde{u}_j(t)=u_j$ and $\tilde{v}_j(t)=v_j$ and set
$\tilde{w}=\sumjn\tilde{u}_j\otimes_C\tilde{v}_j$. Then the function
$t\mapsto\norm{\tilde{w}(t)}$, where $\tilde{w}(t)=\sumjn\tilde{u}_j(t)\otimes
\tilde{v}_j(t)\in X(t)\pro Y(t)$, is not continuous since 
$\norm{\tilde{w}(t_0)}=\norm{w}_{U_0\pro V}>\norm{w}_{U\pro V}=\norm{w(t)}$
if $t\ne t_0$. By last sentence of Theorem \ref{t511} below this discontinuity
implies that $X\pro_CY$ is not normal. (We have used only one direction of
Theorem \ref{t511}, which was deduced in \cite{M6} from a special case in
\cite[Lemma 10]{G}.)
\end{example}

\section{The normal projective tensor product}
 Since for bimodules $X,Y\in\cnrc$ the bimodule $X\pro_CY$ is not necessarily in 
 $\cnrc$, we introduce in this section a new tensor product in the category $\cnrc$. 
 
 We first recall the definition and the construction
 of the normal part of a bimodule.
 \begin{definition}\label{np} Let $A$ be von Neumann algebra.
 The {\em normal part} of a bimodule  $X\in\ara$ is a bimodule 
 $\nor{X}\in\anra$
 together with a contraction $\iota\in\B_A(X,\nor{X})_A$ such that for each
 bimodule $Y\in\anra$ and each $T\in\B_A(X,Y)_A$ there exists a unique
 map $\nor{T}\in\B_A(\nor{X},Y)_A$ such that $\nor{T}\iota=T$ and $\norm{\nor{T}}
 \leq\norm{T}$.
 \end{definition}
 
 By elementary categorical arguments $\nor{X}$ is unique (up to an $A$-bimodule
 isometry) if it exists. To sketch a construction of $\nor{X}$, let 
 $\Phi:A\to\B(\g)$  be the universal representation and
 $\ta=\weakc{\Phi(A)}$  the universal von Neumann
 envelope of $A$. Let $P\in\ta$ be the central projection
 such that the unique weak* continuous extension  of the *-homomorphism 
 $\Phi^{-1}$ has the kernel $\ort{P}\ta$ 
 (see \cite[Section 10.1]{KR} for more details, if necessary).
 Consider $X$ as a subbimodule in its second dual $\bdu{X}$ equipped with
 the canonical bidual $A$-bimodule structure. Since $X$ is representable,
 $\bdu{X}$ can be equipped with a structure of a dual operator $A$-bimodule 
 and by  \cite{B3} or \cite[5.4, 5.7]{BEZ} the bimodule action of $A$ is
 necessarily induced by a pair of $*$-homomorphisms  
 $\pi:A\to A_l(\bdu{X})$ and $\sigma:A\to A_r(\bdu{X})$, where $A_l(\bdu{X})$ and
 $A_r(\bdu{X})$ are certain fixed von Neumann algebras associated to the dual operator
 space $\bdu{X}$ such that $\bdu{X}$ is a normal dual operator $A_l(\bdu{X}),
 A_r(\bdu{X})$-bimodule. Then we may regard $\bdu{X}$ as a normal dual operator
 $\ta$-bimodule through the normal extensions of $\pi$ and $\sigma$ to $\ta$. 
 Now $PXP$ is an $A$-subbimodule
 in $\bdu{X}$, hence so is its norm closure $\nor{X}=\overline{\overline{PXP}}$.
 Finally, define $\iota:X\to \nor{X}$ by $\iota(x)=PxP$. If $T\in\B_A(X,Y)_A$,
 then $\bdu{T}:\bdu{X}\to\bdu{Y}$ is an $\ta$-bimodule map, hence it maps
 $PXP$ into $PYP$. It can be proved \cite{M10} that for a normal bimodule
 $Y\in\anra$ the map 
 $$\iota_Y:Y\to PYP, \ \ \iota_Y(y)=PyP$$
 is isometric, hence we have the factorization $T=\nor{T}\iota_X$,
 where $\nor{T}=\iota_Y^{-1}\bdu{T}|\overline{\overline{PXP}}.$ 
 We summarize the discussion in the following theorem, 
 which is proved in more details in \cite{M10}.

 \begin{theorem}\label{t59}\cite{M10} Let $A$ be a von Neumann algebra, $X\in\ara$
 and regard $X$ as an $A$-subbimodule in $\bdu{X}$. 
 Then  
 $\bdu{X}$ is a normal dual (representable) Banach $\ta$-bimodule
 and  the normal part of $X$ is $\nor{X}=\overline{\overline{PXP}}\subseteq
 \bdu{X}$ with $\iota:X\to\nor{X}$ the map $\iota(x)=PxP$. Moreover,
 \begin{equation}\label{59} \norm{\iota(x)}=
 \inf\left(\sup_j\norm{e_jxf_j}\right), \ \ (x\in X)
 \end{equation} 
 where the infimum is taken over all nets $(e_j)$ and $(f_j)$ of projections
 in $A$ that converge to $1$.
 
 In particular $X\in\anra$ if and only if for all nets of projections $(e_j)$ in
 $A$ and $(f_j)$ in $B$ converging to $1$ we have that $\lim_j\norm{e_jx}=\norm{x}
 =\lim_j\norm{xf_j}$. If $A$ is $\sigma$-finite it suffices to consider increasing 
 sequences of projections instead of nets,
 \end{theorem} 
 We recall that a von Neumann algebra $A$ is  $\sigma$-finite if each orthogonal
 family of nonzero projections in $A$ is countable. 
 The last part of Theorem \ref{t59} was proved for one sided modules in \cite[Theorem 3.3]{M6} and
 this will suffice for our application here since we will consider central 
 bimodules only.

 Now we consider briefly the special case of central bimodules. 
 For a function $f:\Delta\to\br$, let $\essup f$ be the infimum of all
 $c\in\br$ such that the set $\{t\in\Delta:\ f(t)>c\}$ is meager (= contained
 in a countable union of closed sets with empty interiors).
 Define the {\em essential direct sum}, $\eop X(t)$, of a family of Banach spaces 
 $(X(t))_{t\in\Delta}$ as the quotient of the $\ell_{\infty}$-direct
 sum $\oplus_{t\in\Delta}X(t)$ by the zero space of the seminorm 
 $x\mapsto\essup\norm{x(t)}$. Then $\eop X(t)$ with the norm 
 $\dot{x}\mapsto\essup\norm{x(t)}$ is a Banach space and we denote by 
 $e:\oplus_{t\in\Delta}X(t)\to\eop X(t)$ the quotient map. 
 
 \begin{theorem}\cite{M10}\label{t511} Given a bimodule $X\in\coc$ with the canonical decomposition
 $\kappa:X\to\oplus_{t\in\Delta}X(t)$ (see Remark \ref{r47}), its normal part $\nor{X}$ is just the 
 closure of $e\kappa(X)$ in $\eop X(t)$. Moreover, $X\in\cnrc$ if and only
 if for each
 $x\in X$ the function $\Delta\ni t\mapsto\norm{x(t)}$ is continuous.
 \end{theorem}

\begin{definition}\label{d77} If
$X,Y\in\cnrc$, let  $X\npro Y$ be
the completion of $X\otimes_CY$ with the norm
$$\nc(w)=\sup\norm{\phi(w)},\ \ (w\in X\otimes_CY),$$
where the supremum is over all $C$-bilinear contractions $\phi$ from $X\times Y$
into  normal representable $C$-bimodules.
\end{definition}

That $\nc$ is indeed a norm (not just a seminorm) follows since it dominates 
the Haagerup norm on ${\rm MIN}(X)\otimes_C{\rm MIN}(Y)$. (Namely, each completely
contractive bilinear map is contractive. The definiteness
of the Haagerup norm on $X\otimes_CY$ follows from \cite[1.1, 2.3]{M1}). 
We shall
omit the easy proof of the following proposition (the last part of Theorem
\ref{t59} may be used).

\begin{proposition}\label{p78} If $X,Y\in\cnrc$, then $X\npro Y\in\cnrc$
and for each bounded $C$-bilinear map $\psi:X\times Y\to Z\in\cnrc$ there exists a unique $\tilde{\psi}\in{\rm B}_C(X\npro Y,Z)$
such that $\tilde{\psi}(x\otimes_Cy)=\psi(x,y)$ for all $x\in X$, $y\in Y$,
and $\norm{\tilde{\psi}}=\norm{\psi}$.

In particular, $\nc$ is the largest among the norms on $X\otimes_CY$ such that
$\norm{x\otimes_C y}\leq\norm{x}\norm{y}$ and that (the completion of) $X\otimes_CY$
with the norm $\norm{\cdot}$ is a normal representable $C$-bimodule.
\end{proposition}

\begin{proposition}\label{p79} (i) $X\npro Y=\nor{(X\pro_CY)}$, hence the
canonical map 
$$\label{710} X\npro Y\to{\rm ess}\oplus_{t\in\Delta}\left(X(t)\pro 
Y(t)\right)$$
is isometric.

(ii) $\nc(\sumjn x_j\otimes_Cy_j)=\sup\norm{\sumjn\theta(x_j,y_j)},$
where the supremum is over all $C$-bilinear contractions from $X\times Y$ to
$C$.

(iii) $\nc(\sumjn x_j\otimes_Cy_j)=\sup\norm{\sumjn\theta(x_j,y_j)},$
where the supremum is over all  
$\bc$-bilinear $C$-balanced contractions $\theta:X\times Y\to\bc$ such that the map
$C\ni c\mapsto\theta(x,cy)$ is weak* continuous for all $x\in X$, $y\in Y$.
\end{proposition}

\begin{proof} (i) From Proposition \ref{p78} $X\npro Y$ has 
the same universal
property as the normal part of $X\pro_CY$, hence they must be isometric as $C$-bimodules.
Then the rest of (i) follows from Theorem \ref{t511} 
since $(X\pro_CY)(t)=X(t)\pro Y(t)$ by (\ref{72}).

(ii) This is a consequence of the fact that the norm of an element $w$ in
a bimodule $Z\in\cnrc$ is equal to $\sup\{\norm{\phi(w)}:\ \phi\in\B(Z,C),\ 
\norm{\phi}\leq1\}$ (Remark \ref{rmd}).

(iii) For each $w=\sumjn x_j\otimes_Cy_j\in X\otimes_CY$ set
$$\tnc(w)=\sup\norm{\sumjn\theta(x_j,y_j)},$$
where the supremum is over all $\theta$ as in the statement (iii).
Since for each $\rho\in\pd{C}$ of norm $1$ and 
each $C$-bilinear contraction $\phi:X\times Y\to C$ the contraction 
$\theta=\rho\circ\phi$ is $\bc$-bilinear and $C$-balanced,
it follows that $\tnc(w)\geq\nc(w)$. To prove the reverse inequality, regard
a $C$-balanced contraction $\theta:X\times Y\to\bc$ as a linear functional
on $V:=X\pro_CY$.  If the functionals
$c\mapsto\theta_{x,y}(c)=\theta(x,cy)$ are normal, then
$\theta(w)=\lim_j\theta(e_jw)$
for each $w\in X\pro_CY$ and each net of projections $e_j\in C$ converging
to $1$. Thus by (\ref{59}) $\tnc(w)\leq\norm{\iota(w)}$, where $\iota$ is 
the canonical map from $X\pro_CY$ into $\nor{(X\pro_CY)}$. But $\norm{\iota(w)}=
\nu_C(w)$ by (i), hence $\tnc(w)\leq\nc(w)$.
\end{proof}
We call a bimodule $Z\in\cnrc$ {\em strong} if $\sumj p_jz_j
\in Z$ for all bounded sets $(z_j)\subseteq Z$ and orthogonal families of
projections $(p_j)\subseteq C$. (Note that the sum weak* converges in each $\bh$
containing $Z$ as a normal operator $C$-bimodule. Since $Z$ is central, this 
agrees with the definition of general strong bimodules in \cite{M5}.) Strong
modules are characterized as closed in certain topology \cite{M5}, but here we
shall only need  that each bimodule $Z\in\cnrc$ is contained in a
smallest strong bimodule, which follows from \cite[2.2]{M5}.

\begin{remark}\label{r710} Denote by $B_X$ the closed unit ball of a normed space 
$X$. Let $X,Y\in\cnrc$. If $(x_j)_{j\in\bj}\subseteq B_X$, $(y_j)_{j\in\bj}
\subseteq B_Y$ and $(c_j)_{j\subseteq\bj}\subseteq C^+$ are 
such that $\sumj c_j$ weak* converges, then the sum
$\sumj c_jx_j\otimes_Cy_j$
weak* converges in every $\bl$ containing
$X\npro Y$ as a normal $C$-subbimodule since the sum is just the product of bounded operator matrices
\begin{equation}\label{7100}\sumj c_jx_j\otimes_Cy_j=[c_j]^{1/2}_{j\in\bj}{\rm diag}(x_j\otimes_Cy_j)
(c_j^{1/2})_{j\in\bj}.\end{equation}
\end{remark}

\begin{theorem}\label{t711} Given $X,Y\in\cnrc$, let $X\cnpro Y$ be the smallest 
strong $C$-bimodule containing $X\npro Y$. Then every $w\in X\cnpro Y$ 
can be represented in the form
\begin{equation}\label{712} w=\sumj c_jx_j\otimes_Cy_j, \ \  
x_j\in B_X,\ y_j\in B_Y,\ c_j\in C^+, 
\end{equation}
where the sum $\sumj c_j$ weak* converges.
The norm of $w$ is equal to $\inf\norm{\sumj c_j}$ 
over all such representations. 
\end{theorem}

\begin{proof} For $w\in X\npro Y$ 
set $g(w)=\inf\norm{\sumj c_j},$ where the infimum is over all representations
of $w$ as in (\ref{712}). (Since $X\npro Y$ is just the norm completion of 
$X\otimes Y$, a representation of $w$ of the form (\ref{712}) is possible
with the norm convergent series $\sum c_j$.) The inequality
$\nc(w)\leq g(w)$ is proved by essentially the same computation as in the proof 
of Corollary \ref{c73}. The reverse inequality follows from the maximality
of $\nc$ (Proposition \ref{p78}) since  the completion $W$ of 
$X\otimes_CY$  with the norm $g$ is
a representable normal $C$-bimodule. The representability can be verified by
using the characterization of representable bimodules (\cite[Theorem 2.1]{M6},
\cite{Po}) and will be omitted here.
To prove normality 
we may assume that $C$ is $\sigma$-finite for in general
$C$ is a direct sum of $\sigma$ finite algebras and $X$, $Y$ and $X\otimes_CY$ also
decompose into the corresponding direct sums since these are central $C$-bimodules.
If  $W$ is not normal, then by the last part of Theorem
\ref{t59} there exist a sequence of projections $p_j\in C$
increasing to $1$ an element $w\in Z$ and a constant $M$  such that  
$g(p_jw)<M<g(w)$ for all $j$. 
Setting $q_0=p_0$ and $q_j=p_j-
p_{j-1}$ if $j\geq1$, we obtain an orthogonal sequence of projections $q_j$ in $C$
with the sum $1$ such that $g(q_jw)<M$ for all $j$. Thus, for each $j$ we
can choose $x_{ij}\in B_X$ and  $y_{ij}\in B_Y$ and positive elements $c_{ij}\in C$ 
such that  
\begin{equation}\label{pa2} q_jw=\sumi c_{ij}x_{ij}\otimes_Cy_{ij}\ \ 
\mbox{and}\ \ \norm{\sumi c_{ij}}<M,
\end{equation}
where $\bi$ is a sufficiently large index set. Then 
$w=\sum_jq_jw=\sum_j\sumi q_jc_{ij}x_{ij}\otimes_Cy_{ij}$, hence 
(since the projections $q_j$ are central and mutually orthogonal) 
$g(w)\leq\norm{\sum_jq_j\sumi c_{ij}}=\sup_j\norm{\sumi c_{ij}}\leq M$. But this
contradicts the choice of $M$.

To prove that $X\cnpro Y$ consists of elements of the form (\ref{712}),
we may assume (by a direct sum decomposition argument again) that $C$ is $\sigma$ finite.  Then the index set $\bj$ in (\ref{712}) may be taken
to be countable. Given  $w$ as in (\ref{712}),  it follows by 
the Egoroff theorem 
\cite[p. 85]{T} that there exists
an orthogonal sequence of projections $p_k\in C$ with the sum $1$ such that
the sum $\sumj c_jp_k$ is norm convergent for each $k$. Then the sum
$w_k:=\sumj c_jp_kx_j\otimes_Cy_j$ is also norm convergent (to see this, write $w_k$
in the form similar to (\ref{7100})), hence  $w_k\in X\npro Y$ and 
$w=\sum_kw_kp_k\in X\cnpro Y$. Conversely, for each $w\in X\cnpro Y$ there exists an orthogonal 
sequence
of projections $p_k\in C$ such that $wp_k\in X\npro Y$ by \cite[Proposition 2.2]{M5}.
By the first paragraph of the proof $wp_k=\sum_j c_{jk}x_{jk}\otimes_Cy_{jk}$ for some elements
$x_{jk}\in B_X$, $y_{jk}\in B_Y$ and $c_{jk}=c_{jk}p_k\in C^+$ such that 
$\norm{\sum_j c_{jk}}
<\norm{wp_k}+\varepsilon$, where $\varepsilon>0$. 
Then $\norm{\sum_{j,k}c_{jk}}\leq\norm{w}+\varepsilon$ and
$w=\sum_{j,k}c_{jk}x_{jk}\otimes_Cy_{jk}$. This also proves that 
$g(w)\leq\norm{w}$; the reverse 
inequality is clear from (\ref{7100}) by representing $X\cnpro Y$ as
a normal operator $C$-bimodule.
\end{proof}
Since the quotient of a strong bimodule $X\in\cnrc$ by a strong subbimodule
$X_0$ is a strong bimodule in $\cnrc$ by \cite{M10}, we can state the following:

\begin{corollary} If $X_0\subseteq X$ and
$Y_0\subseteq Y$ in $\cnrc$ are strong, then the canonical map
$X\cnpro Y\to(X/X_0)\cnpro(Y/Y_0)$
maps the open unit ball onto the open unit ball.
\end{corollary}

To conclude, we note without presenting the details that results analogous to 
the above ones also hold for the
operator module versions of tensor products (that is, the module versions of
tensor products of operator spaces
studied in \cite{BP}).


\begin{thebibliography}{99}

\bibitem{AP} C. Anantharaman and C. Pop, {\em Relative tensor products and
infinite C$^*$-algebras,} J. Operator Theory {\bf 47} (2002), 389--412.

\bibitem{B3} D. P. Blecher, {\em Multipliers and dual operator algebras,} 
J. Funct. Anal. {\bf 183} (2001), 498--525.
\bibitem{BEZ} D. P. Blecher, E. G. Effros and V. Zarikian, {\em One-sided M-ideals
and multipliers in operator spaces, I,}  Pacific J. Math. {\bf 206} (2002), 287--319. 

\bibitem{BP}  D. P. Blecher and V. I. Paulsen, {\em Tensor products of operator
spaces,}  J. Funct. Anal. {\bf 99} (1991), 262--292.
\bibitem{DF} A. Defant and K. Floret,  {\em Tensor norms and operator ideals,} 
North-Holland Mathematics Studies {\bf 176}, North-Holland Publishing Co., 
Amsterdam, 1993. 

\bibitem{DG} M. J. Dupr\'e and R. M. Gillette, {\em Banach bundles, Banach modules 
and automorphisms of $C\sp{*} $-algebras,} 
Research Notes in Mathematics {\bf 92}, Pitman, Boston, MA, 1983. 

\bibitem{EK}  E.  G. Effros and A.  Kishimoto, {\em Module  maps  and 
Hochschild  -- Johnson cohomology,}  Indiana Univ.  Math.  J. {\bf 36} (1987),
257--276.
\bibitem{ER} E. G. Effros and Z.-J. Ruan, {\em Operator spaces}, London
Math. Soc. Monographs, New Series {\bf 23}, Oxford University Press,
Oxford, 2000.

\bibitem{G} J. Glimm, {\em A Stone-Weierstrass theorem for $C\sp{*} $-algebras,} 
Ann. of Math. {\bf 72} (1960), 216--244.

\bibitem{Ha} H. Halpern, {\em Module homomorphisms of a von Neumann algebra 
into its center}, Trans. Amer. Math. Soc. {\bf 140} (1969), 183--193.

\bibitem{KiR} J. W. Kitchen and D. A. Robbins, {\em Linear algebra in the
category of $C(M)$-locally convex modules,} Rocky Mountain J. Math. {\bf 19}
(1989), 433--480.
\bibitem{KR}  R. V. Kadison and J. R. Ringrose, {\em Fundamentals  of 
the  theory  of operator algebras, Vols. 1, 2},  Academic  Press, 
London, 1983, 1986.
\bibitem{KW} E. Kirchberg and S. Wassermann, {\em Operations on continuous bundles 
of $C\sp *$-algebras,} Math. Ann. {\bf 303} (1995), 677--697.

\bibitem{M1} B. Magajna, {\em The Haagerup norm on the tensor  product 
of operator modules,}  J. Funct. Anal. {\bf 129} (1995), 325--348. 

\bibitem{M5}  B. Magajna, {\em A topology for operator modules over 
W$^*$-algebras,}  J. Funct. Anal. {\bf 154} (1998), 17--41.
\bibitem{M6}  B. Magajna, {\em The minimal operator module of a Banach 
module,} Proc. Edinburgh Math. Soc. {\bf 42} (1999), 191--208.
\bibitem{M10} B. Magajna, {\em Duality and normal parts of operator modules,}
to appear in J. Funct. Anal. {\bf 219} (2005), 306--339. 

\bibitem{P}   V.  I.  Paulsen,  {\em Completely  bounded  maps   and 
operator algebras,} Cambridge Studies in Advanced Mathematics {\bf 78}, Cambridge
University Press, Cambridge, 2002.
\bibitem{Pi} G. Pisier, {\em Introduction to the theory of operator 
spaces}, LMS Lecture Note Series {\bf 294}, Cambridge Univ. Press, Cambridge, 2003.
\bibitem{Po} C. Pop, {\em Bimodules norm\' es repr\' esentables sur des espaces
hilbertiens}, Ph.D. thesis, Universit\' e d'Orl\' eans, January 1999.
\bibitem{R} M. A. Rieffel, {\em Induced Banach representations of Banach algebras 
and locally compact groups} J. Funct. Anal. {\bf 1} (1967), 443--491.
\bibitem{S}  R. R. Smith, {\em Completely bounded module maps  and  the 
Haagerup tensor product,} J. Funct. Anal. {\bf 102} (1991), 156--175.
\bibitem{T}  M. Takesaki, {\em Theory of operator algebras I,} Springer-Verlag,
New-York, 1979.
\bibitem{W1} G. Wittstock, {\em  Extension of completely bounded $C\sp{*} $-module 
homomorphisms,} Operator algebras and group representations, Vol. II (Neptun, 1980), 
238--250, Monogr. Stud. Math., 18, Pitman, Boston, MA, 1984.  



\end{thebibliography}
\end{document}